\documentclass[12pt,twoside,a4paper,reqno]{amsart}
\usepackage{bulletinUPB}
\usepackage{graphics}
\usepackage{amssymb}
\usepackage{eucal}
\usepackage{amsmath}
\usepackage{graphicx}

\newtheorem{preproof}{{\bf Proof }\hspace{-.15cm}}

\newcommand{\qs}{\preceq _s}
\newcommand{\ddd}{\widehat{d}}
\newcommand{\s}{\prec_s}
\newcommand{\ssq}{\subseteq}
\newcommand{\proofbegin}{\begin{preproof}\rm}
\newcommand{\proofend}{\hfill{$\square$} \end{preproof}}
\newtheorem{thm}{Theorem}[section]
\newcommand{\theobegin}{\begin{thm}\rm }
\newcommand{\theoend}{\end{thm}}
\newtheorem{lem}[thm]{Lemma}
\newcommand{\lembegin}{\begin{lem}\rm }
\newcommand{\lemend}{\end{lem}}
\newtheorem{prop}[thm]{Proposition}
\newcommand{\propbegin}{\begin{prop}\rm }
\newcommand{\propend}{\end{prop}}
\newtheorem{conj}[thm]{Conjection. }
\newcommand{\conjbegin}{\begin{conj}\rm }
\newcommand{\conjend}{\end{conj}}
\newtheorem{cor}[thm]{Corollary}
\newcommand{\corbegin}{\begin{cor}\rm }
\newcommand{\corend}{\end{cor}}
\newtheorem{question}[thm]{Question. }
\newcommand{\questionbegin}{\begin{question}\rm }
\newcommand{\questionend}{\end{question}}
\newtheorem{defin}[thm]{definition}
\newcommand{\defbegin}{\begin{defin}\rm }
\newcommand{\defend}{\end{defin}}
\newcommand{\examplebegin}{\begin{preexample}\rm}
\newcommand{\exampleend}{\end{preexample}}

\info{A}{??}{?}{????}
\title{On Maximum Signless Laplacian Estrada Index of Graphs with Given Parameters II}

\author{Ramin Nasiri}
\address{$^1$Department of Mathematics, Faculty of Science, University of Qom, Iran, e-mail: {\tt R.Nasiri.82@gmail.com}}

\author{Hamid Reza Ellahi}
\address{$^2$Department of Mathematics, Faculty of Science, University of Qom, Iran.}

\author{Gholam Hossein Fath-Tabar}
\address{$^3$Department of Mathematics, Faculty of Science, University of Kashan, Iran.}

\author{Ahmad Gholami}
\address{$^4$Department of Mathematics, Faculty of Science, University of Qom, Iran.}

\begin{document}
\pagestyle{headings}
\maketitle

\begin{abstract}


{\it \quad\qu Recently Ayyaswamy \cite{Ayyaswamy01} have introduced a novel concept of the signless Laplacian Estrada index (after here $SLEE$) associated with a graph $G$.
After works, we have identified the unique graph with maximum $SLEE$ with a given parameter such as: number of cut edges,  pendent vertices, (vertex) connectivity and edge connectivity.
In this paper we continue our charachterization for two further parameters; diameter and number of cut vertices.
}
\end{abstract}


\begin{MSC2010}
05C\,12, 05C\,35, 05C\,50.
\end{MSC2010}


\section{Introduction}
Let $G = (V, E)$ be a simple, finite, and undirected graph with vertex set $V(G)$ and the edge set $E(G)$ and $|V(G)|=n$. The adjacency matrix $A=A(G)=[ a_{ij}]$ of $G$ is the binary matrix, where the element $a_{ij}$ is equal to 1 if vertices $i$ and $j$ are adjacent, and 0 otherwise. The matrix $L= D-A$, where $D = diag(d_1,d_2,\ldots,d_n)$ is the diagonal matrix of vertex degrees, is known as the Laplacian matrix of $G$. The matrix $Q=D+A$ is called the signless Laplacian matrix of $G$. We denote the spectrum of $A$, $L$ and $Q$ by $(\lambda_1, \lambda_2,\ldots, \lambda_n)$, $(\mu_1, \mu_2, \ldots, \mu_n)$ and $(q_1, q_2, \ldots, q_n)$, respectively. For a graph $G$, Estrada \cite{est} has defined the \emph{Estrada index} of $G$ as
$$EE(G)=\sum^{n}_{i=1}e^{\lambda_i}.$$
Fath-Tabar et al. \cite{FT} proposed the \emph{Laplacian Estrada index}, in full analogy with Estrada index as
$$LEE(G)=\sum^{n}_{i=1}e^{\mu_i}.$$
Theories of Estrada and  Laplacian Estrada indices of graphs have been extensively studied by several authors (see [2,4-7,9-16,19-26]).

Recently, Ayyaswamy \cite{Ayyaswamy01} developed the innovative notion of the \emph{signless Laplacian Estrada index} as
$$SLEE(G)=\sum^{n}_{i=1}e^{q_i}.$$
He also established lower and upper bounds for $SLEE$ in terms of the number of vertices and edges.
Grone and Merris \cite{Grone01,Grone02} proved that for a bipartite graph $G$, $SLEE(G)=LEE(G)$.

Previousely in \cite{Elahi01}, we characterized the unique graphs with maximum $SLEE$ among the set of all graphs with given number of cut edges,  pendent vertices, (vertex) connectivity and edge connectivity.
In this paper, we continue our research by characterizing the unique graph according to two further parameters: diameter and number of cut vertices.

\section{Preliminaries and Lemmas}
In this section, we first introduce basic definitions, notations and concepts used thorough this paper and restate some proved results found in \cite{Cvetkovic01, Elahi01}.
Then, we prove some needful propositions for proving the main result of the next section.

\defbegin\cite{Cvetkovic01}
A \emph{semi-edge walk} of length $k$ in graph $G$, is an alternating sequence 
$W=v_1 e_1 v_2 e_2 \ldots  v_k e_k v_{k+1}$, where $v_1, v_2, \ldots , v_k, v_{k+1}\in V(G)$, and $e_1, e_2, \ldots , e_k\in E(G)$ such that  the vertices $v_i$ and $v_{i+1}$ are (not necessarily distinct)  end points of edge $e_i$, for any $i=1, 2, \ldots , k$.
 If $v_1=v_{k+1}$, then we say $W$ is a \emph{closed semi-edge walk}.
\defend

By following \cite{Elahi01}, we denote The $k$-th signless Laplacian spectral moment of the graph $G$ by $T_k(G)$ , i.e., $T_k(G)=\sum^{n}_{i=1}q^{k}_{i}$.
\theobegin \cite{Cvetkovic01}\label{theo22}
For a graph $G$, the signless Laplacian spectral moment $T_k$ is equal to the number of closed semi-edge walks of length $k$.
\theoend 
Let $G$ and $G'$ be two graphs, and $x,y\in V(G)$, and $x',y'\in V(G')$.
Denote by $SW_k(G;x,y)$ the set of all semi-edge walks of length $k$ in graph $G$, which are begining at vetex $x$, and ending at vertex $y$. 
For convenience, we use $SW_k(G;x,x)$ instead of $SW_k(G;x)$, and set $SW_k(G)=\bigcup_{x\in V(G)}SW_k (G;x)$.
Thus, by Theorem \ref{theo22}, we have $T_k(G)=|SW_k(G)|$.
Note that, by Taylor expansions, we have
\begin{equation*}\label{eq001}
SLEE(G)=\sum_{k\geq 0}\frac{T_k(G)}{k!}.
\end{equation*}

By $(G;x,y)\qs(G';x',y')$ we mean $ |SW_k (G;x,y)|\leq|SW_k(G'; x', y')|$, for any $k\geq 0$.
Moreover, if  $(G;x,y)\qs(G';x',y')$,
 and there exists some $k_0$ such that $|SW_{k_0}(G;x,y)|<|SW_{k_0} (G';x',y')|$, then we write $(G;x,y)\s(G';x',y')$.
\lembegin\label{lem01}\cite{Elahi01}
Let $G$ be a graph. If an edge $e$ does not belong to $E(G)$, Then $SLEE(G)<SLEE(G+e)$.
\lemend
\lembegin\label{lem02}\cite{Elahi01}
Let $G$ be a graph and $v, u, w_1, w_2, \ldots , w_r\in V(G)$.
suppose that
$E_v=\{e_1=vw_1, \ldots , e_r=vw_r\}$ and 
$E_u=\{e'_1=uw_1, \ldots ,  e'_r=uw_r\}$ are subsets of edges of the complement of $G$.
Let $G_u=G+E_u$ and $G_v=G+E_v$.
If $(G;v)\s (G;u)$, and $(G;w_i,v)\qs (G;w_i,u)$ for each $i=1,2,\ldots,r$,
Then $SLEE(G_v)<SLEE(G_u)$.
\lemend
For a vertex $x$ and an edge $e$, let $SW_k(G; x, [e])$ be the set of all closed semi-edge walks of length $k$ in the graph $G$  starting at vetex $x$ and containing the edge $e$. 
\lembegin\label{lem03}
Let $G$ be a graph and $H=G+e$, such that $e=uv\in E(\overline{G})$. If $(G;v)\qs (G;u)$, then $(H;v)\qs (H;u)$. Moreover, if $(G;v)\s (G;u)$, then $(H;v)\s (H;u)$.
\lemend
\proofbegin
We know that for each $z\in \{u,v\}$, and $k\geq 0$,
$$|SW_k(H;z)|=|SW_k(G;z)|+|SW_k(H;z,[e])|.$$
Since $(G;v)\qs (G;u)$, $|SW_k(G;v)|\leq |SW_k(G;u)|$, for each $k\geq 0$. Thus there is a bijection $f_k:SW_k(G;v)\to A_k\ssq SW_k(G;u)$, for each $k\geq 0$.
\\
It is enough to show that $|SW_k(H;v,[e])|\leq |SW_k(H;u,[e])|$, for each $k\geq 0$. Let $W\in SW_k(H;v,[e])$. We can uniquely decompose $W$ to $W=W_1eW_2e\dots eW_r$, such that $W_i\in SW_{k_i}(G;x,y)$, where $x,y\in \{u,v\}$, and $k_i\geq 0$, and $1\leq i\leq r$. Note that $W_i$ is a semi-edge walk in $G$ and does not contain $e$, Thus the decomposition is unique. For each $W_i$  excactly one of the following cases occurs:
\begin{enumerate}
\item[1)] $W_i\in SW_{k_i}(G;v,v)$. In this case we set $h(W_i)=f_{k_i}(W_i)$. Thus, $h(W_i)\in A_{k_i}\ssq SW_{k_i}(G;u,u)$.
\item[2)] $W_i\in A_{k_i}\ssq SW_{k_i}(G;u,u)$. In this case, set $h(W_i)=f^{-1}_{k_i}(W_i)\in SW_{k_i}(G;v,v)$.
\item[3)] $W_i\in SW_{k_i}(G;u,u)\setminus A_{k_i}$, or $W_i\in SW_{k_i}(G;u,v)$, or $W_i\in SW_{k_i}(G;v,u)$. In these cases, let $h$ fix $W_i$, i.e. $h(W_i)=W_i$.
\end{enumerate}
Now, it is easy to check that the map $h:SW_k(H;v,[e])\to SW_k(H;u,[e])$ by the rule $h_{k}(W)=h_{k}(W_1eW_2e\ldots W_r)=h(W_1)eh(W_2)e\ldots eh(W_r)$ is an injection. 

Note that if there exists $k_0$ such that $|SW_{k_0}(G;v)|<|SW_{k_0}(G;u)|$, then $f_{k_0}$ is not surjective. Thus $h_{k_0}$ is not a surjection, and we have 
$$|SW_{k_0}(H;v,[e])|<|SW_{k_0}(G;u,[e])|$$
which implies that $(H;v)\s(H;u)$.
\proofend

By a similar method, we  prove the following statement:
\lembegin\label{lem04}
Let $G$ be a graph and $H=G+e$, such that $e=uv\in E(\overline{G})$, and $(G;v)\qs (G;u)$.
 If there exists a vertex $x\in V(G)$ such that $(G;x,v)\qs (G;x,u)$, then $(H;x,v)\qs (H;x,u)$.
Moreover, if $(G;v)\s (G;u)$ or $(G;x,v)\s (G;x,u)$, then $(H;x,v)\s (H;x,u)$.
\lemend
\proofbegin
Since $(G;v)\qs (G;u)$, there is a bijection $f_k:SW_k(G;v)\to A_k\ssq SW_k(G;u)$, for each $k\geq 0$. Similarly, since $(G;x,v)\qs (G;x,u)$, there is a bijection $g_k:SW_k(G;x,v)\\ \to B_k\ssq SW_k(G;x,u)$, for each $k\geq 0$.
It is obvious that for each $k\geq 0$,
$$|SW_k(H;x,z)|=|SW_k(G;x,z)|+|SW_k(H;x,z,[e])|$$
where $z\in \{v,u\}$. It is enough to show that for each $k\geq 0$,
$$|SW_k(H;x,v,[e])|\leq |SW_k(H;x,u,[e])|.$$
Let $W\in SW_k(H;x,v,[e])$.
 $W$ decomposes uniquely to $W_1eW_2e\ldots eW_r$, where $W_i$ is a semi-edge walk of length $k_i$ in $G$.
Three cases will be considered as follows for $W_1$:
\begin{enumerate}
\item[1)] If $W_1\in SW_{k_1}(G;x,v)$, then we set $h_1(W_1)=g^{}_{k_1}(W_1)\in B_{k_1}\ssq SW_{k_1}(G;x,u)$.
\item[2)] If $W_1\in B_{k_1}\ssq SW_{k_1}(G;x,u)$, then set $h_1(W_1)=g^{-1}_{k_1}(W_1)\in SW_{k_1}(G;x,v)$.
\item[3)] If $W_1\in SW_{k_1}(G;x,u)\setminus B_{k_1}$, then set $h_1(W_1)=W_1$.
\end{enumerate}
If $1<i\leq r$, then three cases will be considered as follows for $W_i$:
\begin{enumerate}
\item[1)] If $W_i\in SW_{k_i}(G;v)$, then we set $h_i(W_i)=f_{k_i}(W_i)\in A_{k_i}$.
\item[2)] If $W_i\in A_{k_i}\ssq SW_{k_i}(G;u)$, then set $h_i(W_i)=f^{-1}_{k_i}(W_i)\in SW_{k_i}(G;v)$.
\item[3)] If $W_i\in SW_{k_i}(G;u)\setminus A_{k_i}$, or $W_i\in SW_{k_i}(G;v,u)$, or $W_i\in SW_{k_i}(G;u,v)$, then we set $h_i(W_i)=W_i$.
\end{enumerate}
One can check easily that the map $h_k:SW_k(H;x,v,[e])\to SW_k(H;x,u,[e])$ by the rule $h_{k}(W)=h_{k}(W_1eW_2e\ldots W_r)=h_1(W_1)eh_2(W_2)e\ldots eh_r(W_r)$ is injective.

The secound part of the lemma is clear.
\proofend
\section{The graph with maximum SLEE with given diameter}

For $x\in V(G)$, the \emph{eccentricity} $e(x)$ of $x$ is the distance to a vertex of $G$ farthest from $x$, i.e. $e(x)=max\lbrace d(x,y) : y\in V(G)\rbrace$.
The \emph{diameter} $d(G)$ is the maximum eccentricity of the vertices, whereas the
\emph{radius} $r(G)$ is the minimum eccentricity. Also, $x$ is a \emph{central vertex} if
$e(x)=r(G)$ and a \emph{diametral path} is a shortest path between two vertices whose distance is equal to $d(G)$.
For convenience, let us denote $\lceil \frac{d}{2}\rceil$ by $\ddd$ where is the smallest integer number greater than $\frac{d}{2}$.
\\
It is obvious that $K_n$ is the unique graph with diameter $1$.
Also, the path on $n$ vertices $P_n$, is the unique graph with diameter $n-1$.
Furthermore, $K_n-e$ is the graph with maximum $SLEE$ with diameter $2$, where $e$ is an edge of $K_n$. 

\lembegin\label{lem05}
Let $G$ be a graph with diameter $d$, and $P_{d+1}=v_0v_1\ldots v_d$ be a diametral path in $G$. If $d\geq 2$ and $x\in V(G)\setminus V(P_{d+1})$, then $x$ has at most $3$ neighbors in  $V(P_{d+1})$.
\lemend
\proofbegin
Suppose that  $x$ has neighbors $v_{i_1},v_{i_2},\ldots , v_{i_{r}}$ in $P_{d+1}$,  where $r>3$, and $i_1<i_2<\ldots<i_r$. Since $i_r-i_1>2$, the path $P'=v_0v_1\ldots v_{i_1}xv_{i_r}v_{i_r+1}\ldots v_d$ from $v_0$ to $v_d$ is  of  length $d-i_r+i_1+2<d$, which  is a contradiction.
\proofend
Let $n>4$, and $2<d<n-1$, and $1\leq j\leq \ddd$.
 We denote by 
$\mathcal{H}_{d,j}$,
the set of  all graphs $H_{d,j}$, constructed from $K_{n-1-d}$ and $P_{d+1}=v_0v_1\ldots v_d$, by attaching each vertex of $K_{n-d-1}$ to exactly 3 vertices of $P_{d+1}$, such that  for each $x\in V(K_{n-d-1})$, there exists an index $i$, $\ddd-j\leq i\leq \ddd +j-2$, where $x$ is attached to $v_i$, $v_{i+1}$ and $v_{i+2}$.
Therefore, none of $v_i$, $0\leq i<\ddd-j$ or $\ddd +j<i\leq d$, has a neighbor in  $K_{n-d-1}$.
Note that $v_{\ddd}$  is a central vertex of the path $P_{d+1}$. For example, all graphs $H_{4,2}$ with $n=7$  are shown in Fig.1.
\begin{figure}[h]
\center
\includegraphics[width =0.99  \textwidth]{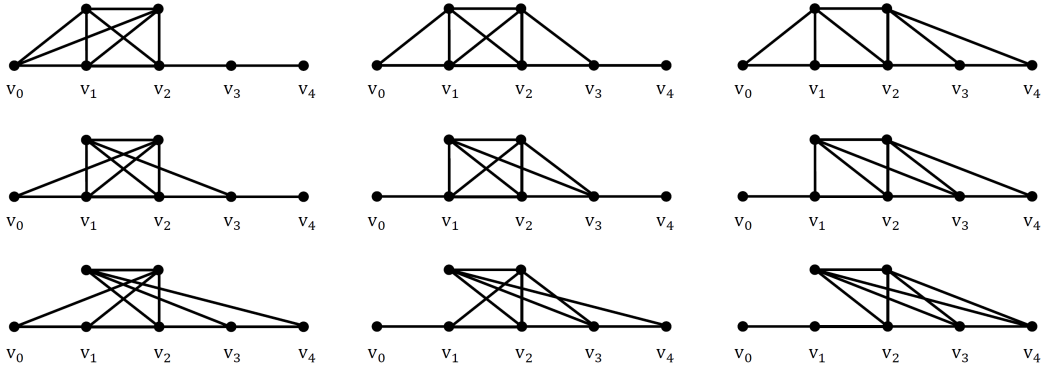}
\caption{All graphs $H_{4,2}$ with $n=7$ .}
\end{figure}
\lembegin\label{lem06}
Let $n>4$, and $2<d<n-1$, and $2\leq j\leq \ddd $. If $H_j\in \mathcal{H}_{d,j}$, then either $H_j\in \mathcal{H}_{d,j-1}$, or there exists a graph, say $H_{j-1}\in \mathcal{H}_{d,j-1}$, such that $SLEE(H_j)<SLEE(H_{j-1})$. 
\lemend
\proofbegin
Let $H_j\in\mathcal{H}_{d,j}$, and $N_K(v_i)=N(v_i)\cap V(K_{n-1-d})$, where $0\leq i\leq d$ and $N(v_i)$ is the set
of vertices that are adjacent to $v_i$.
 To facilitate the understanding of the proof,
  we divide the argument into two parts. We first discuss about  $N_K(v^{}_{\ddd-j})$ and then proceed to  $N_K(v^{}_{\ddd+j})$.  
Note that if $d$ is odd, $j=2$,
 and $N_K(v^{}_{\ddd+2})= \emptyset$, then by renaming the vertices of $P_{d+1}$ such that $v_i$ changes to $v^{}_{d-i}$
  we conclude that $H_j\in \mathcal{H}_{d,j-1}$.
Let $H_j\not\in \mathcal{H}_{d,j-1}$. Therefore either at least one of the vertex subsets $N_K(v^{}_{\ddd-j})$ or $N_K(v^{}_{\ddd+j})$ is not empty, or $d$ is odd and $j=2$ and  $N_K(v^{}_{\ddd+2})$  is not empty.  

If $N_K(v^{}_{\ddd-j})=\emptyset$, then we set $H'_{j-1}=H_j$.
In this case, we have $SLEE(H'_{j-1})=SLEE(H_j)$.
 Let $N_K(v^{}_{\ddd-j})$ be not empty. For convenience, suppose that $v=v^{}_{\ddd-j}$, $y=v^{}_{\ddd-j+1}$, $z=v^{}_{\ddd-j+2}$ and $u=v^{}_{\ddd-j+3}$.
By the definition of $\mathcal{H}_{d,j}$, it is obvious that $N_K(v)\ssq N_K(y)\ssq N_K(z)$, and $N_K(v)\cap N_K(u)=\emptyset$.

Let $E=\{vx:x\in N_K(v)\}$, and $E'=\{ux:x\in N_K(v)\}$, and $H'_j=H_j-E$, and $H'_{j-1}=H'_j+E'$.

By lemma \ref{lem02}, to  show that  $SLEE(H_j)<SLEE(H'_{j-1})$, it is enough to  prove the following statements:
\begin{enumerate}
\item[1)] $(H'_j;v)\s (H'_j;u)$.
\item[2)] $(H'_j;x,v)\qs (H'_j;x,u)$, for each $x\in N_K(v)$.
\end{enumerate}
We start the prove of (1) by the following  claim:\\
{\bf Claim.} $(H'_j;y)\qs (H'_j;z)$:

To prove the claim, let $W\in SW_k(H'_j-e;y)$, where $e=yz$, and $k\geq 0$.
We can decompose $W$ to $W=W_1W_2W_3$, where $W_1$ and $W_3$ are as long as possible and consisting of just the vertices $v_0,v_1,\ldots,y$, and edges in $\{v_tv_{t+1}: 0\leq t\leq \ddd-j\}\cup \{yx:x\in N_K(y)\}$, and $W_2\in SW_{k_2}(H'_j-e;x,w)$,  where $x,w\in N_K(y)\ssq N_K(z)$.
Suppose that $W'_i$ obtains from $W_i$, for $i=1,3$, by replacing each vertex $v_t$ by $v_a$, and  each edge $v_tv_{t+1}$ by $v_av_{a-1}$,  and each edge $yx$ by $zx$, where $x\in N_K(y)$, and $a=2\ddd-2j-t+3$ (In fact, the distance between $v_t$ and $y$ is equal to the distance between $v_a$ and $z$ in $P_{d+1}$).
\\
It is easy to check that the map $f'_k:SW_k(H'_j-e;y)\to SW_k(H'_j-e;z)$ defining by the rule $f'_k(W_1W_2W_3)=W'_1W_2W'_3$ is injective. Thus $(H'_j-e;y)\qs (H'_j-e;z)$. Now, the claim follows from lemma \ref{lem03}.

Let $f_k:SW_k(H'_j;y)\to SW_k(H'_j;z)$ be an injection, for each $k\geq 0$. If $W\in SW_k(H'_j;v)$, then $W$ can be decomposed to $W=W_1W_2W_3$, where $W_2\in SW_{k_2}(H'_j;y)$ is as long as possible . Let $W'_i$ obtain form $W_i$, for each $i=1,3$, by replacing each vertex $v_t$ by $v_a$, and each edge $v_tv_{t+1}$ by $v_av_{a-1}$, where  $a=2\ddd-2j-t+3$.
The map $g_k:SW_k(H'_j;v)\to SW_k(H'_j;u)$, defining by the rule $g_k(W_1W_2W_3)=W'_1f_{k_2}(W_{k_2})W'_3$ is injective.
Note that if $j>2$ or $d$ is even, then the path $v_0v_1\ldots v$ is a proper subgraph of the path $v_dv_{d-1}\ldots u$. Also, if $d$ is odd and $j=2$, then  $N_K(u)\neq \emptyset$, implies that $deg_{H'_j}(v)=2<deg_{H'_j}(u)$. Thus $(H'_j;v)\s (H'_j;u)$ which is (1).

By a similar method used above, we prove the  statement (2). First, we claim that:\\
{\bf Claim.} $(H'_j;x,y)\qs (H'_j;x,z)$, for each $x\in N_K(v)$.

To prove the claim, let $x\in N_K(v)$, and $W\in SW_k(H'_j-e;x,y)$ where $e=yz$. We can decompose $W$ to $W=W_1W_2$ such that $W_1\in SW_{k_1}(H'_j-e;x,w)$ is as long as possible, where $w\in N_K(y)$, and $W_2\in SW_{k_2}(H'_j-e;w,y)$. Suppose that $W'_2$ obtains from $W_2$ by replacing each  vertex $v_t$ by $v_a$, and the edge $wy$ by $wz$, and each edge $v_tv_{t+1}$ by $v_av_{a-1}$, where  $a=2\ddd-2j-t+3$.
\\
One can easily check that the map $h'_k:SW_k(H'_j-e;x,y)\to SW_k(H'_j-e;x,z)$ defining by the rule $h'_k(W_1W_2)=W_1W'_2$ is injective. Thus $(H'_j-e;x,y)\qs (H'_j-e;x,z)$. Now, the claim follows from lemma \ref{lem04}.

Consider $h_k:SW_k(H'_j;x,y)\to SW_k(H'_j;x,z)$ is an injective map, for each $k\geq 0$. Let $W\in SW_k(H'_j;x,v)$. we can decompose $W$ to $W=W_1W_2$, where $W_1\in SW_{k_1}(H'_j;x,y)$ is as long as possible,  and $W_2\in SW_{k_2}(H'_j;y,v)$. Let $W'_2$ obtain from $W_2$ by replacing each vertex $v_t$ by $v_a$, and replacing each edge $v_tv_{t+1}$ by $v_av_{a-1}$, where $a=2\ddd-2j-t+3$.
It is elementary to show that the map $l_k:SW_k(H'_j:x,v)\to SW_k(H'_j;x,u)$ defining by the rule $l_k(W_1W_2)=h_{k_1}(W_1)W'_2$ is an injection.
Thus, $(H'_j;x,v)\qs (H'_j;x,u)$ for each $x\in N_K(v)$. 
It follows the statement (2).
\\
Now, by the above discussion and lemma  \ref{lem02}, we have
$SLEE(H_j)\leq SLEE(H'_{j-1})$,  with equality if and only if $H'_{j-1}=H_j$.
The first part of the argument ends here.

If  $N_K(v^{}_{\ddd+j})$ is empty or $d$ is odd and $j=2$, then $H'_{j-1}\in \mathcal{H}_{d,j-1}$. In this case, set $H_{j-1}=H'_{j-1}$, and of course $SLEE(H_{j-1})=SLEE(H'_{j-1})$.
Let $H'_{j-1}\not\in \mathcal{H}_{d,j-1}$, Then $N_K(v^{}_{\ddd+j})$ is not empty. 
By repeating the above discussion for $v=v^{}_{\ddd+j}$, $y=v^{}_{\ddd+j-1}$, $z=v^{}_{\ddd+j-2}$ and $u=v^{}_{\ddd+j-3}$, we get the graph $H_{j-1}=H'_{j-1}-E+E'$, such that $H_{j-1}\in \mathcal{H}_{d,j-1}$ and  $SLEE(H'_{j-1})<SLEE(H_{j-1})$.
Therefore, $$SLEE(H_j)\leq SLEE(H'_{j-1})\leq SLEE(H_{j-1})\in \mathcal{H}_{d,j-1}$$
with equalities hold, if and only if graphs are equal.
\proofend
The following theorem is our main result of this section, which is determined the  unique graph with maximum $SLEE$ among the set of  all unicyclic graphs with diameter $d$, where $2<d<n-1$.
\theobegin
Let $2<d<n-1$. If $G$ has maximum SLEE with diameter $d$, then
$G=H_{d,1}$.
\theoend
\begin{figure}[h]
\center
\includegraphics[width =0.99  \textwidth]{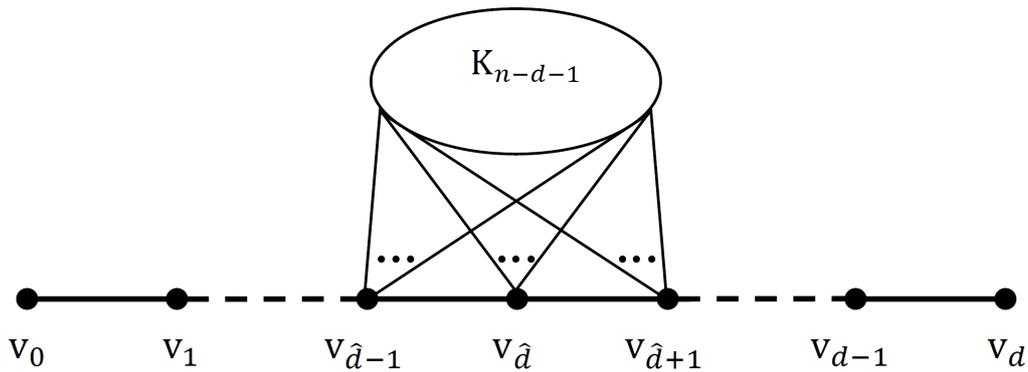}
\caption{The unique graph on $n$ vertices having the maximum SLEE with diameter $d$.}
\end{figure}
\proofbegin
Suppose that $G$ is a graph, having maximum $SLEE$ with diameter $d$.
Let $P_{d+1}=v_0v_1\ldots v_d$ be a diameterical path in $G$, and $H$ be the graph obtained from $G$ by adding some edges such that:
\begin{enumerate}
\item[(a)] $x$ is adjacent with exactly 3 vertices of $P_{d+1}$ in $H$, say $v_i$, $v_{i+1}$ and $v_{i+2}$,  for each $x\in V(G)\setminus V(P_{d+1})$.
\item[(b)] $H-V(P_{d+1})$ is a complete graph on $n-1-d$ vertices.
\end{enumerate}
By lemma \ref{lem05}, such a graph $H$ exists. Obviously, we have $H\in \mathcal{H}_{d,j}$, for some $j$,  $1\leq j\leq \ddd$, and $SLEE(G)\leq SLEE(H)$, with equality if and only if $G=H$.

If $j>1$, then by lemma \ref{lem06}, we may get a sequence of graphs, say $H_{d,j-1}, H_{d,j-2}, \ldots , H_{d,1}$, such that for each $t$, $H_{d,t}\in \mathcal{H}_{d,t}$, and 
\begin{footnotesize}
$$SLEE(G)\leq SLEE(H)\leq SLEE(H_{d,j-1})\leq SLEE(H_{d,j-2})\leq \ldots \leq SLEE(H_{d,1})$$
\end{footnotesize}
with equalities hold, if and only if the graphs are equal.
since the diameter of $H_{d,1}$ is $d$, and $G$  has the maximal $SLEE$ among the set of all graphs with diameter $d$, hence $SLEE(G)=SLEE(H_{d,1})$ which imlplies that $G=H_{d,1}$, as expected.
\proofend

\section{The graph with maximum SLEE with given number of cut vertices}

A \emph{cut vertex} of a graph is a vertex whose removal increases the number of components of the graph. Let $G$ be a connected graph and $x$ be a vertex of $G$, a \emph{block} of $G$ is defined to be a maximal subgraph without cut vertices.\\
A \emph{pendent path} at $x$ in a graph $G$ is a path in which no vertex other than $x$ is incident with any edge
of $G$ outside the path, where $deg_G(x)\geq3$.
In particular, we consider a vertex $x$ as a pendent path at $x$ of length zero in $G$ only when $x$ is neither a pendent vertex nor a cut vertex of $G$.
Let $G$ and $H$ be two vertex-disjoint connected graphs, such that $x \in V (G)$
and $y \in V (H)$. We denote the  \emph{coalescence} of $G$ and $H$ by $G(x) \circ H(y)$, which
is obtained by identifying the vertex $x$ of $G$ with the vertex $y$ of $H$.

\lembegin\label{lem07}
Let $H_1$ and $H_2$ be two graphs and $P_s=y_0y_1\ldots y_{s-1}$ be a path on $s$ vertices,  and $u\in V(H_2)$ and $xy\in E(H_1)$, such that $x\neq y$.
Let $G=\big(H_1(y)\circ P_s(y_0)\big)(x)\circ H_2(u)$. 
If $H_2$ contains a path $Q_{s+2}=ux_1x_2\ldots x_{s+1}$, 
then $SLEE(G)<SLEE(G-E_y+E_{x_1})$,
where $E_y=\big\{yw:w\in N_{H_1}(y)\setminus \{x\}\big\}$ , $E_{x_1}=\big\{x_1w:w\in N_{H_1}(y)\setminus \{x\}\big\}$ and $N_{H_1}(y)$ is the set
of vertices of $H_1$ that are adjacent to $y$ .
\lemend
\begin{figure}[h]
\center
\includegraphics[width =0.99  \textwidth]{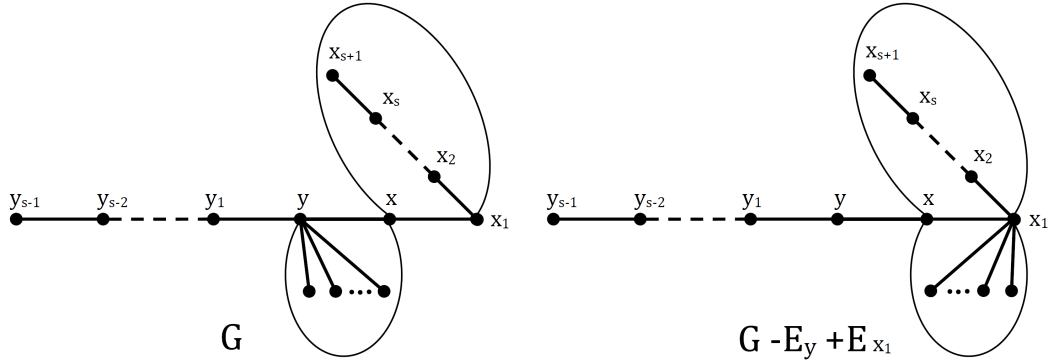}
\caption{An illustration of graphs in Lemma 4.1 .}
\end{figure}
\proofbegin
Let $G'=G-E_y$. By lemma \ref{lem02}, it is enough to show that $(G';y)\s(G';x_1)$, and $(G';w,y)\qs(G';w,x_1)$, for each $w\in  N_{H_1}(y)\setminus \{x\}$.
Let $P'_{s+1}=xy_0y_1\ldots y_{s-1}$, and $A_k=SW_k(G';y)\setminus SW_k(P'_{s+1};y)$, and $B_k=SW_k(G';x_1)\setminus SW_k(Q_{s+2};x_1)$.
Since $P'_{s+1}$ is a proper subgraph of $Q_{s+2}$, it is easy to show that
$|SW_k(p'_{s+1};y)|\leq |SW_k(Q_{s+2};x_1)|$, and for some $k=k_0\geq s$, inequality is strict.
\\
Let $W\in A_k$. 
We may decompose $W$ to $W_1W_2W_3$ such that $W_2\in SW_{k_2}(G';x)$ is as long as possible and $W_1\in SW_{k_1}(G';y,x), W_3\in SW_{k_3}(G';x,y)$ and $k=k_1+k_2+k_3$.
Let $W'_j$ obtain from $W_j$ by replacing each $y_i$ by $x_{i+1}$, where $j=1,3$ and $i=0,1,\ldots , s-1$.
The map $f:A_k\to B_k$ difined by the rule $f(W_1W_2W_3)=W'_1W_2W'_3$ is injective. 
Thus $|A_k|\leq |B_k|$. 
Therefore $|SW_k(G';y)|\leq |SW_k(G';x_1)|$, and for some $k=k_0$ the inequality is strict.
Hence $(G';y)\s (G';x_1)$.

Let $w\in  N_{H_1}(y)\setminus \{x\}$, and $W\in SW_k(G';w,y)$.
We can decompose $W$ uniquely to $W_1W_2$, such that $W_1\in SW_{k_1}(G';w,x)$ is as long as possible.
Let $W'_2$ obtain from $W_2$ by replacing each $y_i$ by $x_{i+1}$, where $W_2\in SW_{k_2}(G';x,y),k=k_1+k_2$ and $i=0,1,\ldots , s-1$.\\
The map $g_{w,k}:SW_k(G';w,y)\to SW_k(G';w,x_1)$ defining by the rule  $g_{w,k} (W_1W_2) = W_1W'_2$ is injective.
Thus $|SW_k(G';w,y)|\leq |SW_k(G';w,x_1)|$ for each $k$.
Therefore $(G';w,y)\qs (G';w,x_1)$, for each $w\in  N_{H_1}(y)\setminus \{x\}$.
\proofend
Let $0 \leq r \leq n-2$. Suppose that $G^r
_n$ is the graph obtained from $K_{n-r}$ by attaching
$n-r$ pendent path of orders $n_1, n_2, \ldots, n_{n-r}$ to its vertices, where
each vertex of $K_{n-r}$ has exactly one pendent path and $\mid n_i-n_j\mid \leq 1$ for
$1\leq i,j\leq n-r$. More precisely, each pendent path is of order $\lfloor \frac{r}{n-r} \rfloor$ or $\lfloor \frac{r}{n-r} \rfloor+ 1$. For example, the graphs $G^r_6$ with $r=0,1,2,3,4$ are shown in Fig.4 .
\begin{figure}[h]
\center
\includegraphics[width =0.99  \textwidth]{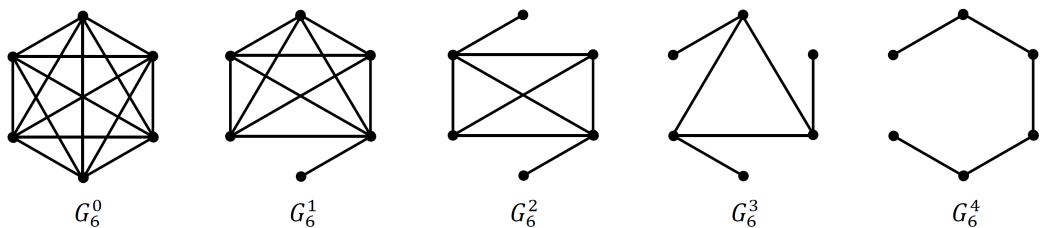}
\caption{The graphs $G^r_6$ with $r=0,1,2,3,4$ .}
\end{figure}
\theobegin
If $0\leq r\leq n-2$, then $G_n^r$ is the unique graph with maximum $SLEE$ among all graphs on $n$ vetices with $r$ cut vertices.
\theoend

\proofbegin
Since $P_n=G_n^{n-2}$ is the uniqe graph with $n-2$ cut vertices, the case $r=n-2$ is obviouse.
If $r=0$, then by lemma \ref{lem01}, $K_n=G_n^0$ is the unique graph on $n$ vertices with maximum $SLEE$.
Let $1\leq r\leq n-3$, and $G$ be a graph with maximum $SLEE$ among all graphs on $n$ vertices with $r$ cut vertices.
\\
First, we prove that $G$ is connected. 
Otherwise, if $G$ is not connected and $x$ is a cut vertex of $G$, then $x$ is also a cut vertex of a component, say $G_1$ of $G$. Let $G_2$ be another component of $G$.
 If $G_2$ has a cut vertex, say $y$, then set $G'=G+\{xy\}$. 
If $G_2$ has no cut vertex, then suppose that $G'$ is the graph obtained from $G$ by attaching $x$ to each vertex of $G_2$.
It is easy to check that in both cases, $G'$ is a graph with $r$ cut vertex and $SLEE(G)<SLEE(G')$, a contradiction. Thus $G$ is connected.

By lemma \ref{lem01}, every block of $G$ is complete.
Let $x$ be a cut vertex contained in at least $3$ blocks, say $B_1$, $B_2$ and $B_3$.
Suppose that, $B_1$ and $B_3$ will be disjointed if the vertex $x$ is removed.
Let $G'$ be the graph obtained from $G$ by attaching each vertex of $B_1$ to each vertex of $B_2$.
Obviously, $G'$ has $r$ cut vertex and by lemma \ref{lem01}, $SLEE(G)<SLEE(G')$, a contradiction.
 Thus, each cut vertex of $G$ is contained in exactly two blocks.
\\
Suppose that $G$ has at least one block with at least 3 vertices.
 Otherwise, since each block of $G$ has 2 vertices, $G$ is a tree with maximum degree 2. 
Thus $G\cong P_n$, and $r=n-2$, a contradiction.
\\
Let $P_s$ be a pendent path with minimum length in $G$ at $x$.
Obviously, $x$ lies in a block of $G$, say $B$, with at least $3$ vertices.
Note that if $s=1$, then $x$ is not a cut vertex.
\\
For each $y\in V(B)$, let $H_y$ be the component of $G-E(B)$ which is containing $y$. 
Obviously, $H_x=P_s$.
Let $y\in V(B)$ such that $y\neq x$.
Let $H$ be the component of $G-\big(E(H_x)\cup E(H_y)\big)$ containing $y$.
We have $G\cong \big(H(x)\circ H_x(x)\big)(y)\circ H_y(y)$.
\\
Suppose that $H_y$ is not a path. 
Since $P_s$ has minimal length, there is a pendent path on at least $s$ vertices at a vertex in $H_y$, say $z$, where $z\neq y$.
Thus $H_y$ contains a path on at least $s+2$ vertices with an end vertex $y$. 
Note that since $H_y$ is not a path, we can choose some vertices of $H_y$ and construct the path of length at least $s+2$ with an end vertex $y$.
By lemma \ref{lem07}, we may get another graph on $n$ vertices with $r$ cut vertices, which has larger $SLEE$, a contradiction.
Therefore,  $H_y$ is a pendent path, say $P_t$ at $y$.
\\
By the choice of $P_s$, we have $t\geq s$.
If $t\geq s+2$, then by lemma \ref{lem07}, we may obtain another graph on $n$ vertices with $r$ cut vertices, which has a larger $SLEE$ that $G$, a contradiction. 
Therefore, for each $y\in V(B)$, $H_y\cong P_s$ or $P_{s+1}$.
Hence $G\cong G_n^r$.
\proofend

\bigskip


\end{document}